\documentclass{amsart}

\usepackage{amssymb}

\usepackage[all]{xy}

\newtheorem{thm}{Theorem}

\newtheorem{lem}[thm]{Lemma}

\newtheorem{prop}[thm]{Proposition}

\theoremstyle{definition}

\newtheorem{say}[thm]{}
\newtheorem{exmp}[thm]{Example}

\newtheorem{ques}[thm]{Question}    

\newtheorem{rem}[thm]{Remark}          

\newtheorem*{ack}{Acknowledgments}      

\newtheorem{defn-thm}[thm]{Definition--Theorem}  
\newtheorem{defn-lem}[thm]{Definition--Lemma}  

\theoremstyle{remark}

\renewcommand{\c}[0]{{\mathbb C}}  

\renewcommand{\o}[0]{{\mathcal O}} 
\newcommand{\z}[0]{{\mathbb Z}}

\renewcommand{\r}[0]{{\mathbb R}} 

\renewcommand{\a}[0]{{\mathbb A}}

\newcommand{\p}[0]{{\mathbb P}}

\newcommand{\q}[0]{{\mathbb Q}}

\newcommand{\qtq}[1]{\quad\mbox{#1}\quad}
\newcommand{\spec}[0]{\operatorname{Spec}}
\newcommand{\pic}[0]{\operatorname{Pic}}

\newcommand{\rank}[0]{\operatorname{rank}}

\newcommand{\discrep}[0]{\operatorname{discrep}}

\newcommand{\codim}[0]{\operatorname{codim}}    
    
\newcommand{\proj}[0]{\operatorname{Proj}}

\newcommand{\cent}[0]{\operatorname{center}}

\newcommand{\sing}[0]{\operatorname{Sing}}

\newcommand{\cl}[0]{\operatorname{Cl}}

\newcommand{\rdown}[1]{\lfloor{#1}\rfloor}

\newcommand{\onto}[0]{\twoheadrightarrow}

\newcommand{\simq}[0]{\sim_{\q}}

\newcommand{\tsum}[0]{\textstyle{\sum}}

\newcommand{\mg}[0]{{\mathbb G}_m}

\def\into{\DOTSB\lhook\joinrel\to}

\def\loccoh#1.#2.#3.#4.{H^{#1}_{#2}(#3,#4)}

\DeclareMathAlphabet{\mathchanc}{OT1}{pzc}%
                                {m}{it}

\begin{document}
\bibliographystyle{amsalpha}


\title{New examples of terminal and 
log canonical singularities}
\author{J\'anos Koll\'ar}

\maketitle


The aim of this note is to revisit the constructions of
\cite{r-c3f, friedman-etal, reid-nndp} to obtain new
examples of terminal and 
log canonical singularities.
First we discuss the general method and then 
work out in detail the following series of examples.

\begin{thm}\label{new.trem.4f.thm} For every $r\geq 4$
 there are germs of 
terminal 4-fold singularities $(0\in X_r)$
such that $K_{X_r}$ is Cartier,
 $X_r\setminus\{0\}$ is simply connected,
 the class group $\cl(X_r)$ is trivial  and
 the embedding dimension of $X_r$ is $r$.
\end{thm}

\begin{thm}\label{new.lc.3f.thm} Let $F$ be a connected 2-manifold 
without boundary. Then there are germs of 
isolated log canonical 3-fold singularities $(0\in X=X(F))$
such that for any resolution  $p:Y\to X$
we have
 $R^1p_*\o_Y\cong H^1(F,\c)$,
 $\pi_1(Y)\cong \pi_1(F)$  and
 $\pi_1(X\setminus \{0\})$ is an extension of $\pi_1(F)$ by a cyclic group.
\end{thm}

\begin{say}[Previous examples]

One of the first results on the singularities of higher dimensional
birational geometry,
proved in \cite{r-c3f}, says  that a 3-fold terminal singularity
whose canonical class is Cartier is a 
hypersurface double point.  Thus every 3-fold terminal singularity
is a quotient of a hypersurface double point by a cyclic group
action. This description can be developed into
a complete list; see \cite{r-ypg}.

For some time it has been an open question if, in some sense,
terminal singularities in higher dimensions also form
an essentially bounded family. The 
first relevant examples are in an unpublished note \cite{bro-rei}.
\medskip

In   earlier examples of germs of isolated log canonical 3-fold singularities
$(0\in X)$ with a resolution  $p:Y\to X$
 we had $\dim R^1p_*\o_Y\leq 2$
(with equality only for cones over an Abelian surface)
and $\pi_1(Y)$  contained a finite index Abelian subgroup.

In  previous higher dimensional  examples we had $\dim R^1p_*\o_Y\leq \dim X-1$
and $\pi_1(Y)$  contained a finite index Abelian subgroup.
These both hold for cones over smooth varieties.

It is also worthwhile to emphasize the difference
between the log canonical and the dlt cases.
Let  $\bigl(0\in X,\Delta\bigr)$ be a germ of a
dlt pair and $p:Y\to X$
any resolution of singularities. Then  
 $R^ip_*\o_Y=0$ for $i>0$ by \cite{elkik} (a simple proof is in
 \cite[5.22]{km-book}) and
  $\pi_1(Y)=1$ by \cite{k-shaf, takayama}.
\end{say}

\subsection*{Singularities with prescribed exceptional divisors}{\ }

Every smooth variety $Z$ can be realized as the exceptional set
of a resolution of an isolated singular point $(0\in X)$; 
one can simply take $X$ to be a cone over $Z$.
More generally, given any scheme $Z$ of dimension $n$, one can ask
if  there is a normal, isolated singularity $0\in X$   of dimension $n+1$
with a  resolution
$$
\begin{array}{ccc}
Z & \subset & Y\\
\downarrow && \downarrow \\
0 & \in & X
\end{array}
\quad \qtq{such that $Y\setminus Z\cong X\setminus \{0\}$.}
$$
An obvious restriction is that $Z$ should have only
hypersurface singularities, but this is not sufficient.
 For instance, consider $Z=(xy=0)\subset \p^3$ and let
$Y$ be any smooth 3-fold containing $Z$. Then the normal bundle
of the line $L:=(x=y=0)\subset Z$ in $Y$ is $\o_{\p^1}(1)+\o_{\p^1}(1)$, thus
$L\subset Y$ deforms in every direction. Hence $Z$ is not  contractible;
it can not even be a subscheme of the exceptional set of a 
resolution of an isolated 3-fold singularity.

It turns out, however, that if we allow $Y$ to be (very mildly)
singular, then  one can construct such $(0\in X)$.

\begin{prop} \label{gen.CI.contains.CI.prop}
Let $P$ be a smooth variety and $Z\subset P$ a subscheme of codimension $n$
that is a local complete intersection. Let $L$ be a line bundle on $P$
such that $L(-Z)$ is generated by global sections.
Let $Z\subset Y\subset P$ be the  complete intersection
of $(n-1)$  general sections of $L(-Z)$. Set
$$
\pi: B_{(-Z)}Y:=\proj_Y\tsum_{m=0}^{\infty}\o_Y(mZ)\to Y.
$$
(We blow up not the ideal sheaf of $Z$ but its inverse in
the  class group.) Then
\begin{enumerate}
\item $B_{(-Z)}Y$ is CM. 
\item $\pi^{-1}_*Z$ is isomorphic to $Z$ and it  is a Cartier divisor
in $ B_{(-Z)}Y$.
\item The exceptional set of $\pi$ is a $\p^1$-bundle of codimension
2 in $B_{(-Z)}Y$.
\item If  $Z$ has only  normal crossing singularities
then $B_{(-Z)}Y$ is terminal  in a neighborhood of $Z$. 
\item If $Z$ has only normal crossing singularities
 and $\dim Z\leq 4$ 
 then
$B_{(-Z)}Y\setminus Z$ is smooth in a neighborhood of $Z$.
\item 
$\omega_{B_{(-Z)}Y}\cong \pi^*\omega_Y\cong 
 \pi^*\bigl(\omega_P\otimes L^{n-1}\bigr)|_Y$ and the 
normal bundle of $\pi^{-1}_*Z$ in $ B_{(-Z)}Y$
is  $\omega_Z\otimes \omega_P^{-1}\otimes L^{1-n}$.
\end{enumerate}
\end{prop}

Proof. The claims (1--5) are \'etale local and, once they are established,
(6) follows from the adjunction formula.
Thus we can assume that $P=\a^N$. Next we write down \'etale local  equations
for $Y$ and for  $B_{(-Z)}Y$ and then read off their properties.

\begin{say}[Local computations]\label{CI.contains.CI.say}
Let $Z\subset \a^N$ be a complete intersection
of codimension $n$  defined
by $f_1=\dots=f_n=0$. Let $Z\subset Y\subset \a^N$ be a general
complete intersection of codimension $n-1$. It is thus defined
by a system of equations
$$
\begin{array}{ccccccc}
h_{1,1} f_1&+ & \cdots & +&h_{1,n} f_n & = & 0\\
 \vdots &&&&\vdots &&\\
h_{n-1,1} f_1&+ & \cdots & +&h_{n-1,n} f_n & = & 0
\end{array}
\eqno{(\ref{CI.contains.CI.say}.1)}
$$
Let $H=(h_{ij})$ be the matrix of the system and 
$H_i$ the submatrix obtained by removing the $i$th column.
Note that for $h_{ij}$ general, the equations
$$
\bigl(\rank H< n-1\bigr)\qtq{that is} \bigl(\det H_1=\cdots=\det H_n=0\bigr)
\eqno{(\ref{CI.contains.CI.say}.2)}
$$
define a codimension 2 subset of $Z$ (\ref{detvar.say}).
If  $f_n({\mathbf x})=0$ then either
$f_1({\mathbf x})=\cdots=f_{n-1}({\mathbf x})=0$ or the system
$$
\begin{array}{ccccccc}
h_{1,1}({\mathbf x})\cdot  y_1&+ & \cdots & +&
h_{1,n-1}({\mathbf x})\cdot   y_{n-1} & = & 0\\
 \vdots &&&&\vdots &&\\
h_{n-1,1}({\mathbf x})\cdot   y_1&+ & \cdots & +&
h_{n-1,n-1}({\mathbf x})\cdot   y_{n-1} & = & 0
\end{array}
\eqno{(\ref{CI.contains.CI.say}.3)}
$$
has a nontrivial solution $y_i=f_i({\mathbf x})$,
 thus $\det H_n({\mathbf x})=0$. We can do this for any
$j$ instead of $n$, hence we get that
$$
(f_j=0)=Z\cup (f_j=\det H_j=0)\subset Y.
\eqno{(\ref{CI.contains.CI.say}.4.j)}
$$
By our argument, this holds set-theoretically, but since
$Y$ is CM, (\ref{CI.contains.CI.say}.4.j) in fact holds scheme-theoretically.
The formula (\ref{CI.contains.CI.say}.4.j) also suggests that computing
$$
\pi: B_{(-Z)}Y:=\proj_Y\tsum_{m=0}^{\infty}\o_Y(mZ)\to Y
$$
is the same as blowing up the ideal $(f_j,\det H_j)$ for any $j$:
$$
B_{(-Z)}Y\cong B_{(f_j,\det H_j)}Y.
\eqno{(\ref{CI.contains.CI.say}.5.j)}
$$
Again, because of the possible difference between the powers of
the ideal  $(f_j,\det H_j)$ and its symbolic powers, 
so far  we only know that
the $S_2$-hull of $B_{(f_j,\det H_j)}Y$ is $B_{(-Z)}Y$.
Next we  show that $B_{(f_j,\det H_j)}Y$ is CM
hence $S_2$, thus
(\ref{CI.contains.CI.say}.5.j) indeed holds.

$Z$ is CM, hence (\ref{CI.contains.CI.say}.4.j)
and (\ref{liaison.lem}) imply that $(f_j=\det H_j=0)$ is CM
and so is $(f_j=\det H_j=0)\times \p^1$.
In $Y \times \p^1$ the equation $\bigl(sf_j=t\det H_j\bigr)$
defines  $B_{(f_j,\det H_j)}Y \cup \bigl((f_j=\det H_j=0)\times \p^1\bigr)$.
Again using (\ref{liaison.lem}) we see that $B_{(f_j,\det H_j)}Y$ is CM.

The formula (\ref{CI.contains.CI.say}.5.j) shows that
$\pi$ is an isomorphism whenever $f_j\neq 0$ or
$\det H_j\neq 0$. Letting $j$ vary, the first set of these
conditions define $Z$ and the second set 
defines a codimension 2 subset of $Z$ by (\ref{CI.contains.CI.say}.2).

For notational simplicity, let us compute the $j=n$ case.
In $\a^N\times \p^1_{st}$ the blow-up satisfies the equations
(\ref{CI.contains.CI.say}.1) and  $sf_n=t\det H_n$. 
Multiplying the system (\ref{CI.contains.CI.say}.1) by the
determinant-theoretic adjoint of $H_n$, we get the equations
$$
\det H_n\cdot  \bigl(f_1,\dots, f_{n-1}\bigr)^{tr}  + 
 f_n\cdot H_n^{adj}\cdot  \bigl(h_{1,n},\dots, h_{n-1,n}\bigr)^{tr}
 =  0.
\eqno{(\ref{CI.contains.CI.say}.6)}
$$
 Multiplying through by $s$,
substituting $sf_n=t\det H_n$ and dividing by $\det H_n$
we get new equations for the blow-up:
$$
s\cdot  \bigl(f_1,\dots, f_{n-1}\bigr)^{tr} +
t\cdot H_n^{adj}\cdot  \bigl(h_{1,n},\dots, h_{n-1,n}\bigr)^{tr}
 =  0
\eqno{(\ref{CI.contains.CI.say}.7)}
$$
These, together with   $sf_n=t\det H_n$
show that $(t=0)\subset B_{(f_n, \det H_n)}Y$
is isomorphic to $Z$. 
Furthermore, $H_n^{adj}\cdot  \bigl(h_{1,n},\dots, h_{n-1,n}\bigr)^{tr}$
is the zero vector exactly where $\rank H<n-1$, proving (3).

Thus $U:=(s\neq 0)\subset B_{(f_n, \det H_n)}Y$ is an open
 neighborhood of $Z$.  In $U$,   the equations
(\ref{CI.contains.CI.say}.1) are consequences of
(\ref{CI.contains.CI.say}.7) and  $sf_n=t\det H_n$.
Thus the $n$ equations 
(\ref{CI.contains.CI.say}.7) and  $sf_n=t\det H_n$
define $U$, hence it is a
local complete intersection.

If $Z$ has hypersurface singularities, then we can set
$f_1, \dots, f_{n-1}$ to be linear and independent.
Thus the equations (\ref{CI.contains.CI.say}.7) can be used
to eliminate variables. Possibly after shrinking
$U$ and choosing new \'etale coordinates,   we end up with one equation
$$
U:=(f_n=t\det H_n)\subset \a^{N-n+2}_{({\mathbf x},t)},
\eqno{(\ref{CI.contains.CI.say}.8)}
$$
where $H_n$ is a general $(n-1)\times (n-1)$ matrix
of polynomials in the ${\mathbf x}$-variables.

If $Z$ has  normal crossing singularities then
$Z$ is slc, but if $E$ is a  divisor over $Z$
with discrepancy $\leq 0$ then $\cent_ZE$ is either
a stratum of $Z$ (if $\discrep(E,Z)=-1$)
or is a codimension 1 point in a stratum of $Z$ (if $\discrep(E,Z)=0$)
(cf.\ \cite[2.29]{km-book}).
Thus, by the precise inversion of adjunction \cite[17.3, 17.12]{k-etal},
$U$ is terminal near $Z$ iff it is terminal at the 
codimension 1 points of the strata of $Z$.
At these points, the equation  (\ref{CI.contains.CI.say}.8)
is
$$
\bigl(x_1\cdots x_m=t\bigr)
\qtq{or}
\bigl(x_1\cdots x_m=tx_{m+1}\bigr)\subset \a^{N-n+2}_{({\mathbf x},t)};
\eqno{(\ref{CI.contains.CI.say}.9)}
$$
these are terminal singularities. 
For later purposes we also note the following.
\medskip

{\it Claim \ref{CI.contains.CI.say}.10. Assume that  $Z$ has  normal crossing 
singularities and $L$ is ample. 
\begin{enumerate}
\item [a)] Every irreducible component
of $\sing Z$ contains a point where, in suitable local analytic coordinates,
$$
\bigl[Z\subset B_{(-Z)}Y\bigr] \cong \bigl[(t=0)\subset 
(x_1 x_2=tx_{3})\bigr].
$$
\item [b)] $B_{(-Z)}Y$ is smooth  at a general point of every stratum of $Z$.
\qed
\end{enumerate}
}
\medskip

If $\dim Z\leq 4$ then $\dim \sing Z\leq 3$ and
by (\ref{detvar.say}) we may assume that every point
of $Z\subset B_{(-Z)}Y$ is described by a local equation
$\bigl(x_1\cdots x_m=t\bigr)$ or $\bigl(x_1\cdots x_m=tx_{m+1}\bigr)$
 for some $m\leq 4$. This shows (5).

Starting with $\dim Z= 5$, we get  singularities outside $Z$
of the form
$$
\bigl(x_1x_2=t(x_3x_4-x_5x_6)\bigr)\subset \a^7_{({\mathbf x},t)}.
$$
This is still of type $cA$. If $\dim Z= 6$ then we also get
triple points of the form
$$
\bigl(x_1x_2x_3=t(x_4x_5-x_6x_7)\bigr)\subset \a^8_{({\mathbf x},t)}. \qed
$$
\end{say}

\begin{say}[Determinantal varieties]\label{detvar.say}
We have used the following properties of determinantal varieties;
see \cite{bru-vet} for a general treatment.

Let $V$ be a smooth, affine variety,
$V_i\subset V$ a finite set of  smooth, affine subvarieties
and ${\mathcal L}\subset \o_X$
a finite dimensional base point free linear system.
Let  $H_{n,m}=\bigl(h_{ij}\bigr)$ be an $n\times m$ matrix
whose entries are general elements in ${\mathcal L}$. Then
 for every $i$
\begin{enumerate}
\item the singular set of $V_i\cap \bigl(\det H_{n,n}=0\bigr)$
 has codimension 4 in $V_i$ and
\item the set $V_i\cap \bigl(\rank H_{n,n-1}<n-1\bigr)$
 has codimension 2 in $V_i$.
\end{enumerate}
\end{say}

The following is a basic observation of liaison theory 
(cf.\ \cite[Sec.21.10]{eis-ca}).

\begin{lem}\label{liaison.lem} Let $W$ be a Gorenstein scheme
of pure dimension $n$ 
that is a union of two of its closed 
subschemes $W_1, W_2$ of pure dimension $n$. 
If $W_1$ is CM then so is $W_2$.
\end{lem}

Proof. Set $D:=W_1\cap W_2$. There is an exact sequence
$$
0\to \o_{W_2}(-D)\to \o_{W}\to \o_{W_1}\to 0.
$$
Here $\o_{W_2}(-D)$ is the largest subsheaf of $\o_{W}$
whose support is in $W_2$.
This shows that $\o_{W_2}(-D)$ is CM. 
Similarly, we have
$$
0\to \o_{W_1}(-D)\to \o_{W}\to \o_{W_2}\to 0.
$$
Hom it to $\omega_W$ to get
$$
0\to \omega_{W_2}\to \omega_{W}\to \omega_{W_1}(D)\to 0.
$$
Here again $\omega_{W_2}$ is the largest subsheaf of $\omega_{W}$
whose support is in $W_2$. Since $W$ is  Gorenstein,
$\o_W\cong \omega_{W}$ which implies that
$\omega_{W_2}\cong \o_{W_2}(-D)$. Thus $\omega_{W_2}$ is CM and
so is $\o_{W_2}$. \qed
\medskip

Let now $Z$ be a projective, connected,  local complete intersection scheme
of pure dimension $n$ and $L$ an ample line bundle on $Z$.
A large multiple of $L$ embeds $Z$ into $P:=\p^N$ for some $N$.
Applying (\ref{CI.contains.CI.say}) we get $Z\subset B_{(-Z)}Y$ where the
normal bundle of $Z$ is very negative. Thus $Z$ can be contracted
(analytically or as an algebraic space) and we get
a singular point $(0\in X)$. Its properties are
summarized next.

\begin{thm} \label{CI.sch.is.exc.set.thm}
Let $Z$ be a projective, connected,  local complete intersection scheme
of pure dimension $n$ and $L$ an ample line bundle on $Z$.
Then for $m\gg 1$ there are 
  germs of  normal  singularities  $\bigl(0\in X=X(Z,L,m)\bigr)$
with a partial resolution
$$
\begin{array}{ccc}
Z & \subset & Y\\
\downarrow && \hphantom{\pi}\downarrow\pi \\
0 & \in & X
\end{array}
\quad\qtq{where $Y\setminus Z\cong X\setminus\{0\}$}
$$
such that
\begin{enumerate}
\item  $Z$   is a Cartier divisor
in $Y$,
\item the normal bundle of $Z$ in  $Y$
is  $\omega_Z\otimes L^{-m}$,
\item if $Z$ is snc then 
 $Y$ has only terminal  singularities
and (\ref{CI.contains.CI.say}.10.a--b) hold.
\item if  $\dim Z\leq 4$ then $(0\in X)$ is an isolated singular point. \qed
\end{enumerate}
\end{thm}

Next we consider various properties of these $(0\in X)$. 

\begin{prop} \label{almost.cones.properties.1}
Let   $0\in X$ be  a  germ of a normal  singularity
with a partial resolution
$$
\begin{array}{ccc}
Z & \subset & Y\\
\downarrow && \hphantom{\pi}\downarrow\pi \\
0 & \in & X
\end{array}
\quad\qtq{where $Y\setminus Z\cong X\setminus\{0\}$.}
$$
Assume that   $Z$ is a reduced, slc, Cartier divisor
with normal bundle $L^{-1}$ where $L$ is ample and $-K_Z$ is nef. Then:
\begin{enumerate}
\item $R^i\pi_*\o_Y\cong H^i(Z,\o_Z)$.
\item The restriction $\pic(Y)\to \pic(Z)$ is an isomorphism.
\item If $Y$ is smooth then 
$\cl(X)= \pic(Y)/\bigl\langle [Z_i]:i\in I\bigr\rangle$
where the $\{Z_i:i\in I\}$ are the irreducible components of $Z$.
\item If $\o_Z(-K_Z)\cong L^r$ then
\begin{enumerate}
\item $K_X$ is Cartier and $(0\in X)$ is lc,
\item If $r>0$ and $Y$ has canonical singularities then
$(0\in X)$ is canonical.
\item If $r>1$, $Y$ has canonical singularities
and $Y\setminus Z$ has terminal  singularities then
$(0\in X)$ is terminal.
\end{enumerate}
\end{enumerate}
\end{prop}

Proof. We  think of a germ as a small analytic neighborhood
of $0\in X$. To be precise, we assume that
$X$ is Stein and $Z$ is a deformation retract of $Y$.
This in particular implies that we have isomorphisms
$$
H^*\bigl(Z, \z\bigr)\cong H^*\bigl(Y, \z\bigr)\qtq{and}
\pi_1(Z)\cong \pi_1(Y).
\eqno{(\ref{almost.cones.properties.1}.5)}
$$
Let $I_Z\subset \o_Y$ be the ideal sheaf of $Z\subset Y$.
Then the completion of $R^i\pi_*\o_Y$ equals the inverse limit
of the groups $H^i(Y,\o_Y/I_Z^m)$, hence it is enough to prove that
the maps
$$
H^i\bigl(Y,\o_Y/I_Z^{m+1}\bigr)\to H^i\bigl(Y,\o_Y/I_Z^{m}\bigr)\to 
\cdots \to H^i(Z,\o_Z)
$$ 
are all isomorphism. For each of these we have 
an exact sequence
$$
 H^i\bigl(Z,L^{m}\bigr)\to 
H^i\bigl(Y,\o_Y/I_Z^{m+1}\bigr)\to H^i\bigl(Y,\o_Y/I_Z^{m}\bigr)\to
H^{i+1}\bigl(Z,L^{m}\bigr)
$$
and $H^i\bigl(Z,L^{m}\bigr)=0$ by Kodaira vanishing, proving (1).
(We need  Kodaira vanishing in a slightly more general
setting that usual  \cite[2.4]{km-book}. 
The normal crossing case is easy to derive by induction, or see
\cite[9.12, 12.10]{shaf-book} for the general log canonical setting
using \cite{k-db}.)

A similar argument proves (2) using the exact sequence
$$
0\to L^m\to \bigl(\o_Y/I_Z^{m+1}\bigr)^*\to \bigl(\o_Y/I_Z^{m}\bigr)^*\to 1.
$$
which gives 
$$
H^1\bigl(Z,L^{m}\bigr)\to
\pic\bigl(\spec_Y\o_Y/I_Z^{m+1}\bigr)\to\pic\bigl(\spec_Y\o_Y/I_Z^{m}\bigr).
$$
If  $Y$ is smooth then $\cl(Y)= \pic(Y)$ which implies (3).

Finally, by adjunction, $\bigl(K_Y-(r-1)Z\bigr)|_Z\sim 0$,
hence  $K_Y-(r-1)Z\sim 0$ by (2). 
Thus $K_X\sim \pi_*\bigl(K_Y-(r-1)Z\bigr)$ is Cartier,
$K_Y\sim \pi^*K_X+(r-1)Z$ and the rest of (4) is clear.
\qed

\begin{prop} \label{almost.cones.properties.2} 
Let   $Z$ be a connected, reduced, proper,  Cartier divisor in a 
normal analytic space
$Y$ such that $Z$ is a deformation retract of $Y$. 
 For each irreducible component $Z_i\subset Z$
let $\gamma_i\subset Y\setminus Z$ be a small loop around $Z_i$.
\begin{enumerate}
\item Assume that $\codim_Z (Z\cap \sing Y)\geq 2$. Then 
the $\gamma_i$ (and their conjugates) generate 
$\ker\bigl[\pi_1\bigl(Y\setminus Z\bigr)\to \pi_1(Y)\bigr]$.
\item  Assume in addition that $Z$ is a normal crossing scheme and
\begin{enumerate}
\item every codimension 1 stratum  of $Z$ contains a node of $Y$ and
\item every codimension 2 stratum  of $Z$ contains a smooth point 
of $Y$.
\end{enumerate}

\noindent Then $\ker\bigl[\pi_1\bigl(Y\setminus Z\bigr)\to \pi_1(Y)\bigr]$
is cyclic and is generated by any of the $\gamma_i$.
\end{enumerate}
\end{prop}

Proof. Let
$\rho^0:(Y\setminus Z)^{\sim}\to {Y\setminus Z}$
be an \'etale cover that is trivial on all the $\gamma_i$.
This means that $\rho^0$ extends to an
 \'etale cover of $Y\setminus (Z\cap \sing Y)$. 
By a Lefschetz type theorem \cite[XIII.2.1]{sga2}
it then extends to an \'etale cover of $Y$.
This proves (1).

Next let 
$\rho^0:(Y\setminus Z)^{\sim}\to {Y\setminus Z}$
be an \'etale cover  and $\rho:Y^{\sim}\to Y$ its extension
to a normal, ramified cover of $Y$.
Assume that there is at least one point $y\in \rho^{-1}(Z)$
where $\rho$ is \'etale. We claim that $\rho$ is then everywhere \'etale.

First we show that $\rho$ is \'etale
generically on every irreducible component of $\rho^{-1}(Z)$
and on every codimension 1 stratum  of $\rho^{-1}(Z)$.
Since $\rho^{-1}(Z)$ is connected in codimension 1,
we only need to show that if $\tilde Z_i, \tilde Z_j$ are
two irreducible components of $\rho^{-1}(Z)$ such that
 they intersect in codimension 1 and $\rho$ is 
generically  \'etale along  $\tilde Z_i$ then it is also
generically  \'etale along  $\tilde Z_j$ and
along $\tilde Z_i\cap \tilde Z_j$.

Here we use the existence of nodes (2.a). 
In local coordinates we have
$$
(Z\subset Y)\cong (t=0)\subset (x_1x_2=tx_3)\subset \a^{n+1}_{(x_1,\dots, x_n,t)}.
$$
Note that $(x_1x_2=tx_3)\setminus (t=0)\sim \c^*\times \c^{n-1}$,
thus, in this neighborhood, the two loops
$\gamma_1$ around  $(t=x_1=0)$ and $\gamma_2$ around  $(t=x_2=0)$ are homotopic.
Thus any \'etale cover of $(x_1x_2=tx_3)\setminus (t=0)$
that is unramified along $(t=x_1=0)$ is also
unramified along $(t=x_2=0)$.

At a  general point of a codimension 2 stratum  of $Z$ 
we use (2.b). Thus, in  
local coordinates we have
$$
(Z\subset Y)\cong (t=0)\subset (x_1x_2x_3=t)\subset 
\a^{n+1}_{(x_1,\dots, x_n,t)}
$$
and a cover that is \'etale along $Z$ is trivial.
Finally, since $Z$ is a hypersurface singularity, 
it is locally simply connected at codimension $\geq 3$ points
\cite[X.3.4]{sga2}.
Thus any (possibly ramified) cover of $Z$ that is
\'etale  outside a subset of codimension $\geq 3$
is everywhere \'etale.
 \qed

\begin{rem} Note that the seemingly artificial condition
(\ref{almost.cones.properties.2}.2.a) is necessary,
even if $Y$ is smooth everywhere.

As an example, let $S$ be a resolution of a rational surface
singularity $(0\in T)$ with exceptional curve $E\subset S$. Take
$Z:=E\times \p^1\into S\times \p^1=: Y$. 

Then $Z$ is simply connected but
$\pi_1(Y\setminus Z)\cong \pi_1(S\setminus E)$
is infinite, nonabelian  as soon as $T$ is not a quotient
singularity.

The presence of nodes along the double
locus of $Z$  (\ref{gen.CI.contains.CI.prop}.10),
 which at first seemed to be a  blemish of the
construction,
thus turned out to be of great advantage to us.
\end{rem}

\begin{say} \label{getting.rid.of.kernel.say}
 The kernel of
$\bigl[\pi_1\bigl(Y\setminus Z\bigr)\to \pi_1(Y)\bigr]$
in (\ref{almost.cones.properties.2})
can be infinite cyclic; for instance this happens
if  $Z$ is an abelian variety.

However, if  $\pi_1(Z)=1$ and $Z\subset Y$ is contractible then
 $\pi_1\bigl(Y\setminus Z\bigr)$ is finite.
To see this note that since  $\pi_1\bigl(Y\setminus Z\bigr)$ is
abelian, it is enough to show that $\gamma_i$ has finite order in
$H_1\bigl(Y\setminus Z, \z\bigr)$.
We can now repeatedly cut $Y$ by hyperplanes until
it becomes a smooth surface, hence a resolution of a normal
surface singularity. The latter  case is
 computed in \cite{mumf-top}.

An especially simple situation is when
$\omega_Z\cong L^r$ for some $r$.
We can choose $L$ to be non-divisible in $\pic(Z)$.
Then the normal bundle of $Z$ in $Y$ is
$L^{r-m}$, thus $\o_Y(-Z)$ is divisible by
$(m-r)$ in 
$\pic(Y)/\bigl\langle [Z_i]:i\in I\bigr\rangle$.
We can thus take a degree $(m-r)$ cyclic cover of
$X\setminus \{0\}$ and replace $Z\subset Y$ with another
diagram
$$
\begin{array}{ccc}
Z & \subset & \tilde Y\\
\downarrow && \hphantom{\pi}\downarrow\pi \\
0 & \in &\tilde X
\end{array}
\quad
$$
where $\tilde Y\setminus Z\cong \tilde X\setminus\{0\}$,
the normal bundle of $Z$ is  $L^{-1}$ and
$ \tilde X\setminus\{0\}$ is simply connected.

The singularities of $\tilde Y$ are, however, a little worse
than before. 
The argument after (\ref{CI.contains.CI.say}.8)
shows that they are canonical, not terminal.
At double points of $Z$, the  original local equations
$(x_1x_2=t)$ or $(x_1x_2=tx_3)$ become
$$
(x_1x_2=s^{m-r})\qtq{or} (x_1x_2=s^{m-r}x_r).
$$
\end{say}

\subsection*{Construction of 
log canonical 3-fold singularities}{\ }

By (\ref{CI.sch.is.exc.set.thm}--\ref{almost.cones.properties.2}), 
the following construction
of reducible snc surfaces implies (\ref{new.lc.3f.thm}). 
The construction and the proof of its properties
were clearly known to the authors of \cite{friedman-etal},
though only some of it is there explicitly.

\begin{prop}\label{snc.CY.surf.prop}
 For every connected 2-manifold without boundary $F$ 
there are (many) 
connected algebraic surfaces $Z=Z(F)$  with snc singularities such that
\begin{enumerate}
\item  $K_Z\sim 0$  if $F$ is orientable and  $2K_Z\sim 0$ if $F$ is not 
orientable,
\item $h^i(Z,\o_Z)= h^i(F,\c)$ and
\item $\pi_1(Z)\cong \pi_1(F)$.
\end{enumerate}
\end{prop}

We first describe the irreducible components of these surfaces
and then explain how to glue them together.

\begin{say}[Rational surfaces with an anticanonical cycle] \label{r.s.a.c}
Fix $m\geq 1$ and 
let $Z$ be a rational surface such that $-K_{Z}$ is linearly equivalent
to a length $m$
cycle of rational curves $C_1,\dots, C_{m}$.
One can get such surfaces by starting with 3 lines in $\p^2$,
then repeatedly blowing up an intersection point and adding the exceptional 
curve to the collection of curves $\{C_j\}$.

Let $L$ be an ample line bundle on $Z$ and consider the sequence
$$
0\to L(-\tsum C_j)\to L\to 
L|_{\sum C_j}\to 0.
$$
Since $-\sum C_j\sim K_Z$, we see that 
$H^1\bigl(Z, L(-\tsum C_j)\bigr)=0$.
Thus we have a surjection
$$
H^0\bigl(Z, L\bigr)\onto 
H^0\bigl(\tsum C_j, L|_{\sum C_j}\bigr).
$$
Since $\pic^0(\sum C_j)\cong \mg$, $L|_{\sum C_j}$ has a section
which has exactly one zero on each $C_j$ (not counting multiplicities).
Thus there is a divisor $A_L\in |L|$ such that 
 $A_L\cap C_j$ is a single point for every $j$. 

We would like to choose $L$ such that $\deg L_{C_j}$ is independent of $j$.
This is not always possible, but it can be arranged a follows.

Assume that  the self intersections $\bigl(C_j\cdot C_j\bigr)$
are $\leq -2$ for every $j$.
(This can be achieved by blowing up points on the $C_j$ if necessary.)
Then, for any $j$, all the other curves
form the resolution of a cyclic quotient singularity,
and their intersection form is negative definite.
Thus if $H_j$ is any ample divisor on $Z$ then there is an effective
linear combination
$$
H'_j:=H_j+\tsum_{i\neq j} a_iC_i
$$
such that $\bigl(H'_j\cdot C_i)=0$ for $i\neq j$ and 
$\bigl(H'_j\cdot C_j\bigr)>0$.
Set
$$
H:=\tsum_j \frac1{(H'_j\cdot C_j)}\cdot H'_j.
$$
$H$ is 
 an ample $\q$-divisor  that has degree 1 on every $C_j$.
A suitable multiple gives the required line bundle $L$.

Let next $P\subset \r^2$ be a convex polygon with vertices
$v_1,\dots, v_m$ and sides $S_i=[v_i, v_{i+1}]$.
We map $P$ into the algebraic surface $Z$ as follows.

We map the vertex $v_i$ to the point $C_i\cap C_{i-1}$.
We can think of $C_i\cong \c\p^1$ as a sphere with 
$C_i\cap C_{i+1}$ as north pole, $C_i\cap C_{i-1}$ as south pole and
the unique point $A_L\cap C_i$ as a point on the equator.
We map the side $S_i$ to a semicircle in $C_i\cong \c\p^1\sim S^2$
whose midpoint is $A_L\cap C_i$. Since $\pi_1(Z)=1$, this mapping of the
boundary of $P$ extends to   $\tau:P\to Z(\c)$ (whose image could have
self-intersections). 
\end{say}

\begin{say}[Gluing rational surfaces with  anticanonical cycles]
 \label{r.s.a.c.g}
Let $F$ be a (connected)
topological surface and $T$ a triangulation of $F$.

Dual to $T$ is a subdivision of $F$ into polygons
$P_i$ such that at most 3 polygons meet at any  point.

For each  polygon $P_i$ with  sides $S_1^i,\dots, S_{m_i}^i$
(in this cyclic order) choose  
 a rational surface
$Z_i$ with an anticanonical 
cycle of rational curves $C_1^i,\dots, C_{m_i}^i$ (in this cyclic order)
and a map  $P_i\into Z_i(\c)$ as in (\ref{r.s.a.c}). 

If the sides $S_j^i$ and $S_k^l$ are identified on $F$
by an  isometry $\phi_{jk}^{il}:S_j^i\to S_k^l$,
then there is a unique isomorphism $\Phi_{jk}^{il}:C_j^i\to C_k^l$
extending $\phi_{jk}^{il}$.

These gluing data define a surface $Z=\cup_i Z_i$ and the maps
$\tau_i$ glue to 
$\tau:F\to Z(\c)$. 
Since only 3 polygons meet at any  vertex, we get an snc surface.
The curves $H_i$
glue to an ample Cartier divisor $H$ on $Z$.

We claim that $\tau$ induces an isomorphism
$\pi_1(F)\to \pi_1(Z(\c))$.
This is clear for the 1-skeleton where each 1-cell in $F$
is replaced by a $\c\p^1\sim S^2$. 

As we attach each polygon $P_i$ to the 1-skeleton,
we kill an element of the fundamental group corresponding to
its boundary. On each rational surface $Z_i$ with anticanonical cycle
$\sum_j C^i_j$ we have a surjection
$$
\pi_2\bigl(Z_i, \tsum_j C^i_j\bigr)\onto \pi_1\bigl(\tsum_j C^i_j\bigr),
$$
thus as we attach  $Z_i$ we kill the same element of the fundamental group.
Thus $\pi_1(F)\cong \pi_1(Z(\c))$.

The statement about the homology groups is proved in 
\cite[pp.26--27]{friedman-etal}. 
\end{say}

\medskip


\subsection*{Construction of 
terminal 4-fold singularities}{\ }

Let $W$
be a smooth Fano 3-fold of index 2. That is,
$W$ is smooth and there is an ample line bundle $L$
such that $-K_W\sim L^2$. Then the cone  
$$
C(W,L):=\spec \tsum_{m\geq 0} H^0\bigl(W, L^m\bigr)
$$
 is a terminal singularity.

Such  smooth Fano 3-folds have been classified, they
give examples only up to  embedding dimension $7$.

As a generalization, one can try to look for
 Fano 3-folds of index 2 with  normal crossing singularities.
Thus the irreducible components of its normalization
are normal crossing pairs  $(W_i, S_i)$ with
an ample divisor $H_i$
such that $-\bigl(K_{W_i}+S_i\bigr)\sim 2H_i$.

The first examples of such pairs that come to mind are 
$\bigl(\p^3, \p^1\times \p^1\bigr)$ and
$\bigl(\q^3, \p^1\times \p^1\bigr)$.
These are all the examples where the underlying variety $W$ is also Fano.

In (\ref{indxe.2.logafno.exmp}) we construct an
infinite sequence of   index 2 log Fano pairs
$\bigl(P_r, \p^1\times \p^1\bigr)$ where
$P_r$ is a $\p^2$-bundle over $\p^1$. 

Any 2 of these examples can be glued together
to get infinitely many families of 
 Fano 3-folds of index 2 with  normal crossing singularities.

Then we  apply (\ref{gen.CI.contains.CI.prop})
and (\ref{CI.sch.is.exc.set.thm}--\ref{almost.cones.properties.2})
 to conclude the proof of (\ref{new.trem.4f.thm}).

\begin{exmp}\label{indxe.2.logafno.exmp} For  $r\geq 0$ set
$$
\pi: P_r:=\proj_{\p^1}\bigl(E_r\bigr)\to \p^1
\qtq{where} E_r:=\o_{\p^1}+\o_{\p^1}+\o_{\p^1}(r).
$$
Write $\o_P(a,b):=\o_{P_r}(a)\times \pi^*\o_{\p^1}(b)$.
Since $\pi_*\o_P(a,b)=S^aE_r\otimes \o_{\p^1}(b)$, we see that
$\o_P(a,b)$ is ample iff $a>0$ and $b>0$. 
Let $S_r\subset P_r$ be the surface corresponding to the unique section
of $\o_P(1,-r)$. 
Note that  $K_{P_r}\sim \o_P(-3, r-2)$, 
$S_r\cong \proj_{\p^1}\bigl(\o_{\p^1}+\o_{\p^1}\bigr)\cong
\p^1\times \p^1$ and
$$
-\bigl(K_{P_r}+S_r\bigr)\sim \o_P(2,2)\qtq{is ample on $P_r$. }
$$ 
For any $a,b\geq 0$ there are natural isomorphisms
$$
H^0\bigl(P_r, \o_P(a,b)\bigr)=
H^0\bigl(\p^1, S^aE_r\otimes\o_{\p^1}(b)\bigr)
$$
and $S^aE_r\otimes\o_{\p^1}(b)$ naturally decomposes
as the sum of line bundles of the form
$\o_{\p^1}(cr+b)$ where $0\leq c\leq a$. This makes it  easy to compute
the spaces $H^0\bigl(P_r, \o_P(a,b)\bigr)$
 and to show the following
\begin{enumerate}
\item For $a, b\geq 0$ 
the restriction maps
$$
H^0\bigl(P_r, \o_P(a,b)\bigr)\to H^0\bigl(S_r, \o_S(a,b)\bigr)
\qtq{are surjective. }
$$
\item For $a_i, b_i\geq 0$ 
the multiplication maps $$
H^0\bigl(P_r, \o_P(a_1,b_1)\bigr)\otimes H^0\bigl(P_r, \o_P(a_2,b_2)\bigr)
\to H^0\bigl(P_r, \o_P(a_1+a_2,b_1+b_2)\bigr)
$$
are surjective. 
\end{enumerate}
\end{exmp}

\begin{say}[Fano 3-folds of index 2]\label{index.2.sing.fanos}
We consider in detail two series of examples.

1. Fix  $r\geq 0$ and let
$Z_{r}$ be obtained from
$\bigl(P_r, S_r\bigr)$ and $\bigl(\p^3, S_0\bigr)$
by an isomorphism $S_r\cong \p^1\times \p^1\cong S_0$.
Then
$\omega_{Z_{r}}$ is ample and isomorphic to  $L^{-2}_{r}$ where
$L_r$ is ample.

One can check that  $ h^0(Z_{r}, L_{r}\bigr)=r+6$ and 
 the algebra
$\sum_{m\geq 0} H^0(Z_{r}, L^m_{r}\bigr)$
is generated by $ H^0(Z_{r}, L_{r}\bigr)$.

2. Fix  $r,s\geq 0$ and let
$Z_{rs}$ be obtained from
$\bigl(P_r, S_r\bigr)$ and $\bigl(P_s, S_s\bigr)$
by an isomorphism $S_r\cong \p^1\times \p^1\cong S_s$.
Then
$\omega_{Z_{rs}}$ is ample and  isomorphic to $L^{-2}_{rs}$ where
$L_{rs}$ is ample.
Using (\ref{indxe.2.logafno.exmp}.1--2) we compute that
 $ h^0(Z_{rs}, L_{rs}\bigr)=r+s+8$ and 
\ the algebra
$\sum_{m\geq 0} H^0(Z_{rs}, L^m_{rs}\bigr)$
is generated by $ H^0(Z_{rs}, L_{rs}\bigr)$.

For  $r\in \{0,1\}$, the $(Z_r,L_r)$ series should
give degenerations of smooth Fano 3-folds.
The simplest is $r=0$.
Take $\p^1\times \p^2$ and embed it into $\p^5$ by
$\o(1,1)$. Under this embedding, 
$\p^1\times \p^1\subset \p^1\times \p^2$ becomes a quadric;
its is thus contained in a 3-plane $H^3$. 
The union of $\p^1\times \p^2$ and of the 3-plane 
gives $Z_0$. It is a $(2,2)$ complete intersection.

The construction of $X_r$ can also be realized by
taking the cone over $Z_r$ and then deforming the (reducible) cone by 
taking high-order perturbations of the two quadratic
defining equations.
\end{say}

Putting these together with 
(\ref{CI.sch.is.exc.set.thm}) and 
(\ref{getting.rid.of.kernel.say}) we get the following.

\begin{prop} \label{terminal.4f.new.prop}
Let $(Z,L)$ be any of the pairs from (\ref{index.2.sing.fanos}).
Then there are 4-dimensional, normal,
isolated  singularities $0\in X$ 
with a partial resolution
$\pi:(Z\subset Y)\to (0\in X)$
such that
\begin{enumerate}
\item  $Y$ has only canonical singularities of type $cA$,
\item $Z\subset Y$ is a Cartier divisor and its
normal bundle is $L^{-1}$,
\item $0\in X$ is terminal, $K_{X}$ is Cartier, 
$X\setminus \{0\}$ is simply connected and
\item the embedding dimension of $(0\in X)$ is $h^0(Z,L)$.\qed
\end{enumerate}
\end{prop}

It remains to compute the class group of $ X$.
We apply (\ref{almost.cones.properties.1}.3), which needs a resolution of
singularities of $Y$. Our set-up is simple enough
that this can be done explicitly.

\begin{say}[Explicit resolution]\label{Explicit.resolution.say}
Let $Y$ be  a 4-fold and $Z\subset Y$ a Cartier divisor.
Assume that $Z$ is the union of 2 smooth
components $Z=Z_1\cup Z_2$ and along the intersection
$S:=Z_1\cap Z_2$ in suitable local analytic coordinates 
$\bigl[Z_1\cup Z_2\subset Y\bigr]$
is isomorphic to 
$$
\begin{array}{lll}
\bigl[(x_1=s=0)\cup (x_1=s=0) \subset (x_1x_2=s^m)\bigr]&\subset& \a^5
\qtq{or}\\
\bigl[(x_1=s=0)\cup (x_1=s=0) \subset (x_1x_2=s^mx_3)\bigr]&\subset& \a^5.
\end{array}
$$
Let $C\subset S$ be the curve defined locally by
$(x_1=x_2=x_3=s=0)$. 

We can resolve the singularities  by iterating the
following steps.
\begin{enumerate}
\item  If $m\geq 3$, we blow up $S$. We get 2 exceptional divisors,
both are $\p^1$-bundles over $S$ with 2 disjoint sections. In the above local
coordinates the equation changes to $(x_1x_2=s^{m-2})$ or $(x_1x_2=s^{m-2}x_3)$.
\item   If $m=2$, we blow up $S$; the resulting 4-fold is smooth.
 We get 1 exceptional divisor,
which is a conic bundle over $S$ with 2 disjoint sections. It is isomorphic to a
 $\p^1$-bundle over $S$ blown up along $C$ (contained in one of the sections).
\item  If  $m=1$,  we blow up 
the component, call it $N$,  that
intersects $Z_1$.  
(This component is $Z_2$ iff there are no previous blow-ups.)
The birational transform
of $Z_1$ is isomorphic to $Z_1$ and the  birational transform
of $N$ is isomorphic to $N$ blown-up along $C$.
\end{enumerate}

Thus at the end we have a chain of smooth 3-folds
$$
E_0:=Z_1, E_1, \dots, E_{m-1}, E_m:=Z_2
$$
such that
\begin{enumerate}\setcounter{enumi}{3}
\item the only intersections are $S_i:=E_i\cap E_{i+1}\cong S$ for
$0\leq i<m$.
\item $E_2,\dots, E_{m-1}$ are $\p^1$-bundles over $S$ with 2 disjoint sections
\item $E_1$ is a $\p^1$-bundle over $S$  blown up along $C$ 
if $m\geq 2$ and $Z_2$  blown up along $C$ 
if $m=1$.
\end{enumerate}
By taking the cohomology of the exact sequence
$$
0\to \z_{\cup_i E_i}\to 
\tsum_{i=0}^m \z_{E_i}\to \tsum_{i=0}^{m-1} \z_{S_i} \to 0
$$
we get an exact sequence
$$
\tsum_{i=0}^{m-1}   H^1\bigl(S_i, \z\bigr)\to
H^2\bigl(\cup_i E_i, \z\bigr)\to 
\tsum_{i=0}^m   H^2\bigl(E_i, \z\bigr)\to
\tsum_{i=0}^{m-1}   H^2\bigl(S_i, \z\bigr).
$$
Assume now that  
$$
H^1\bigl(S, \z\bigr)=0\qtq{and}
H^2\bigl(Z_2, \z\bigr)\to H^2\bigl(S, \z\bigr)
\qtq{is surjective.}
\eqno{(\ref{Explicit.resolution.say}.7)}
$$
 Setting $h^i(\ ):=\rank H^i(\ , \z)$  we get that
$$
\begin{array}{rcl}
h^2\bigl(\cup_i E_i\bigr)&=&
\tsum_{i=0}^m   h^2\bigl(E_i\bigr)-
\tsum_{i=0}^{m-1}   h^2\bigl(S_i\bigr)\\[1ex]
&=& h^2\bigl(Z_1\bigr)+
h^2\bigl(Z_2\bigr)-h^2\bigl(S\bigr)+h^2\bigl(C\bigr)+(m-1).
\end{array}
$$

In (\ref{index.2.sing.fanos}) $S\cong \p^1\times \p^1$ and 
$C$ is connected,
 thus the formula becomes
$$
h^2\bigl(\cup_i E_i\bigr)=h^2\bigl(Z_1\bigr)+h^2\bigl(Z_2\bigr)+(m-2).
$$
Therefore, by (\ref{almost.cones.properties.1}.3), the class group of $X$,
obtained by contracting $Z\subset Y$ to a point,  satisfies
$$
\begin{array}{rcl}
\rank \cl(X)&\leq& h^2\bigl(Z_1\bigr)+h^2\bigl(Z_2\bigr)+(m-2)-(m+1)\\
&=&h^2\bigl(Z_1\bigr)+h^2\bigl(Z_2\bigr)-3.
\end{array}
\eqno{(\ref{Explicit.resolution.say}.8)}
$$
Thus we see that for the series  $(0\in X_r)$
we get $\cl(X_r)=0$ but for the series  $(0\in X_{rs})$
we get $\cl(X_{rs})\cong\z$.
\end{say}

\begin{rem} The series $X_{rs}$ can also be constructed as follows.
Set
$$
Y_{rs}:=\proj_{\p^1}\bigl(\o_{\p^1}+\o_{\p^1}+\o_{\p^1}(r)+\o_{\p^1}(s)\bigr).
$$
This visibly contains both $P_r, P_s$ as divisors and
$K_{Y_{rs}}+P_r+ P_s\sim \o_{Y_{rs}}(-2,-2)$.
We can now take the affine cone  $C_a\bigl(Y_{rs},\o_{Y_{rs}}(1,1) \bigr)$ and
inside it the cone  $C_a\bigl(P_r+ P_s \bigr)$. 
Perturbing the equation of the  cone  $C_a\bigl(P_r+ P_s \bigr)$
as in \cite[2.43]{km-book}
we get our varieties $X_{rs}$.

Since  $\cl(X_{rs})=\z$, the singularities
$X_{rs}$  have a small modification; now we can see this
explicitly. The cone $C_a\bigl(Y_{rs}\bigr)$
has a small resolution determined by the pencil
of the 4-planes that are the cones over the fibers of $Y_{rs}\to \p^1$.

By contrast, I believe that 
one can not realize the series $Z_r$ as hypersurfaces
in a smooth variety.
\end{rem}

\begin{rem} Although these examples show that 
terminal 4-fold singularities do not form an
``essentially bounded family'' in 
the most naive sense, it should not be considered a
final answer.

There are several ways to ``simplify'' a given 
terminal  singularity  $(0\in X)$.

First, if $X\setminus\{0\}$ is not simply connected,
one can pass to the universal cover; see (\ref{more.exmp.ques.5}).
Since finite group actions on a given singularity are frequently
not too hard to understand, it makes sense to concentrate
on those  singularities  $(0\in X)$ for which
$X\setminus\{0\}$ is  simply connected.

A harder to use reduction step is the following.
If there is a Weil divisor $D\subset X$ that is not
$\q$-Cartier, then there is a (unique) small  modification
$g:X_D\to X$ such that $g^{-1}_*D$ is 
$\q$-Cartier and $g$-ample. This $X_D$ also has terminal singularities.
Even for 3-folds this reduction is quite subtle  since the
existence of such a $D$ depends on the
coefficients of the equations.
(For instance, if $X_f:=(x_1x_2+f(x_3,x_4)=0)\subset \a^4$
then  $\cl(X_f)=0$  iff $f(x_3,x_4)$ is irreducible (as a power series);
see \cite[2.2.7]{k-etc}.)

Nonetheless this suggests that the
most basic terminal  singularities  $(0\in X)$
of dimension $n$ are those that satisfy both
$\pi_1\bigl(X\setminus\{0\}\bigr)=1$ and
$\cl(X)=0$.

Finally,  one might argue that any collection of examples
constructed in a simple uniform way forms an
``essentially bounded family.''
 In essence, I have just replaced the old
recipe  ``take a cone over a smooth Fano 3-fold of index 2''
with a newer  recipe  ``take a cone over an snc Fano 3-fold of index 2
and deform it.''
The complete list of
snc Fano 3-folds of index 2  (and also many higher dimensional cases)
is in  \cite{fujita-11}.
Thus my examples do form an ``essentially bounded family.''
\end{rem}

\subsection*{Questions}{\ }

These examples raise several questions.

\begin{ques} \label{more.exmp.ques.1}
Is there a more conceptual way to construct
all these examples?
\end{ques}

Comments: A quite general approach could be the following.

Let $0\in X\subset \a^N$ be a subscheme and assume
that the deformations of all singularities of $X\setminus\{0\}$
are unobstructed.  
(For instance, $X$ could be a cone over a projective variety
whose singularities have unobstructed deformations.)
Then all the obstructions to deforming
$X$ are supported at the origin, hence the obstruction space
for $X$ should be finite dimensional.
On the other hand, if $0\in X$ is not an isolated singularity,
then the deformation space itself should be infinite dimensional.

We can thus expect that deformations of $X$, even those deformations
that induce a flat deformation of the tangent cone at $0\in X$,
give pretty much a ``general'' deformation at the singularities
of $X\setminus\{0\}$.

A technical difficulty in carrying this out is that
infinite dimensional deformation spaces are difficult to
handle. It is also not clear what to expect over $X\setminus\{0\}$.

As an example, assume that $X\subset \a^N$
given by a determinantal condition $\bigl(\rank H_{r,r+1}<r\bigr)$
as in (\ref{detvar.say}).
Then every deformation of $X$  is also determinantal.
The singular set is then given by the condition
$\bigl(\rank H_{r,r+1}<r-1\bigr)$, which is a subset of codimension 6.
We aim to keep the origin still singular, thus, if $\dim X\geq 7$,
we can not smooth $X\setminus\{0\}$.

As we saw in (\ref{gen.CI.contains.CI.prop}),
 we can still expect for a general deformation to have
very nice singularities outside the origin.
I do not know any definite general results.

\begin{ques}\label{more.exmp.ques.2} Let $Z$ be a projective snc scheme with
$K_Z\sim 0$. Is there any restriction
on $\pi_1(Z)$?
\end{ques}

Comments:  We can try to follow the construction on (\ref{r.s.a.c.g})
but already in dimension 3 this seems quite difficult.

Assume that we have a polyhedron $P$ with boundary $F\sim S^2$.
As in (\ref{r.s.a.c.g}), we can use $F$ and its
triangulation to build a surface $Z$ such that $K_Z\sim 0$.
The next step would be to find a 3-fold $X$ that contains $Z$
as a divisor such that $K_X+Z\sim 0$. It is not clear that this is 
always possible.
Even if for each polyhedron $P$ such a 3-fold $X(P)$ exists,
gluing them together is probably quite a bit more subtle
than for surfaces.

Here is a quite interesting simple case.

Let $P$ be the dodecahedron. For each face we choose a degree 5
del Pezzo surface; it contains a length 5 chain of $-1$-curves
that form an anticanonical cycle. We can glue these together
to get a surface $Z$ with 12 irreducible components
such that $K_Z\sim 0$. I do not know how to construct a
rationally connected 3-fold $X$  containing $Z$
as a divisor such that $K_X+Z\sim 0$.

This example would be quite interesting since many hyperbolic 3-manifolds
admit a tiling with copies of the dodecahedron.

\medskip

By an observation of \cite{simp}, for each finitely presented group
$\Gamma$ there is a seminormal scheme $Z(\Gamma)$ 
such that  $\pi_1\bigl(Z(\Gamma)\bigr)\cong \Gamma$.

 The following variant was explained  to me by  M.~Kapovich.

Start with a finite simplicial complex $C$ whose fundamental group is $\Gamma$.
For each $k$-simplex $c \in C$ we take
$Z(c):=\c\p^k$. Use the incidence relation in $C$ to glue the spaces
$Z(c)$ together (using linear embeddings). The result is a singular
projective scheme $Z(\Gamma)$ whose fundamental group is $\Gamma$.

If $C$ is a topological manifold, then $ K_{Z(\Gamma)}\sim 0$ and
the singularities of $Z(\Gamma)$ are
simple normal crossing in codimension $1$ but 
more complicated in codimension $\geq 2$.
In codimension $2$ we get degenerate cusps, and their deformation
theory is quite subtle \cite{ghk}. 
I do not know if these examples can be realized as 
exceptional sets of partial resolutions of log canonical singularities.

Kapovich also outlined an argument of how to modify this
construction to obtain a 
projective variety $Y(\Gamma)$ with simple normal crossing singularities
whose fundamental group is $\Gamma$. However, the canonical
class of $Y(\Gamma)$ is not trivial.

\begin{ques} \label{more.exmp.ques.3}
Let $(0\in X)$ be the germ of a log canonical
singularity and $g:Y\to X$ a resolution.
 Is there any restriction
on $\pi_1(Y)$?
\end{ques}


\begin{ques} \label{more.exmp.ques.4}
Let $(X,D)$ be an lc pair with $K_X+D\simq 0$.
What are the possible groups $\pi_1(D)$?
Is $\pi_1(D)$ a birational invariant of the pair $(X,D)$?
What about $\pi_1\bigl(\rdown{\Delta}\bigr)$
if  $(X,\Delta)$ is lc and  $K_X+\Delta\simq 0$?
\end{ques}

Comments: Let $(X,D)$ be an lc pair with $K_X+D\simq 0$.
If $D\neq 0$ then $X$ is uniruled and if $D$ is rationally chain connected
then, by looking at the MRC fibration we conclude that $X$ is
rationally  connected. This implies that $H^i(X, \o_X)=0$ for $0<i$.
By looking at the exact sequence
$$
0\to \omega_X\cong \o_X(-D)\to \o_X\to \o_D\to 0
$$
we conclude that $H^i(D, \o_D)=0$ for $0<i<\dim D$.
It is natural to hope that $\pi_1(D)$ should be finite.

\begin{ques} \label{more.exmp.ques.5}
Let $(0\in X)$ be the germ of a log terminal singularity.
Is $\pi_1(X\setminus \{0\})$ finite?
\end{ques}

Comments: This is probably the most basic of the above questions.
It is closely related to the following old problem:

Let $(X, \Delta)$ be a dlt Fano variety. Is
the fundamental group of the smooth locus $X^{ns}$ finite?
The answer is yes in dimension 2 but even that case needs work
\cite{gur-zha, fkl, kee-mck}.

 \begin{ack} This paper was written while I visited  RIMS, Kyoto University.
I thank  S.~Mori and S.~Mukai for the invitation and their hospitality.
I am grateful to K.~Fujita,  M.~Kapovich, S.~Mori, C.~Simpson  and Ch.~Xu for
many   comments and corrections.
Partial financial support  was provided by  the NSF under grant number 
DMS-0758275.
\end{ack}


\begin{thebibliography}{GHK11}

\bibitem[BR10]{bro-rei}
Gavin Brown and Miles Reid, \emph{Anyone knows these guys?},
  http://dl.dropbox.com/u/10909533/anyoneknows.pdf, 2010.

\bibitem[BV88]{bru-vet}
Winfried Bruns and Udo Vetter, \emph{Determinantal rings}, Lecture Notes in
  Mathematics, vol. 1327, Springer-Verlag, Berlin, 1988. \MR{953963
  (89i:13001)}

\bibitem[Eis95]{eis-ca}
David Eisenbud, \emph{Commutative algebra}, Graduate Texts in Mathematics, vol.
  150, Springer-Verlag, New York, 1995, With a view toward algebraic geometry.
  \MR{1322960 (97a:13001)}

\bibitem[Elk81]{elkik}
R.~Elkik, \emph{Rationalit{\'{e}} des singularit{\'{e}}s canoniques}, Inv.
  Math. \textbf{64} (1981), 1--6.

\bibitem[FKL93]{fkl}
Akira Fujiki, Ryoichi Kobayashi, and Steven Lu, \emph{On the fundamental group
  of certain open normal surfaces}, Saitama Math. J. \textbf{11} (1993),
  15--20. \MR{1259272 (94m:32042)}

\bibitem[FM83]{friedman-etal}
Robert Friedman and David~R. Morrison (eds.), \emph{The birational geometry of
  degenerations}, Progr. Math., vol.~29, Birkh\"auser Boston, Mass., 1983.
  \MR{690262 (84g:14032)}

\bibitem[Fuj11]{fujita-11}
Kento Fujita, \emph{(in preparation)}, 2011.

\bibitem[GHK11]{ghk}
M.~{Gross}, P.~{Hacking}, and S.~{Keel}, \emph{{Mirror symmetry for log
  Calabi-Yau surfaces I}}, arXiv1106.4977G, 2011.

\bibitem[Gro68]{sga2}
Alexander Grothendieck, \emph{Cohomologie locale des faisceaux coh\'erents et
  th\'eor\`emes de {L}efschetz locaux et globaux {$(SGA$} {$2)$}},
  North-Holland Publishing Co., Amsterdam, 1968, Augment\'e d'un expos\'e par
  Mich\`ele Raynaud, S\'eminaire de G\'eom\'etrie Alg\'ebrique du Bois-Marie,
  1962, Advanced Studies in Pure Mathematics, Vol. 2. \MR{0476737 (57 \#16294)}

\bibitem[GZ95]{gur-zha}
R.~V. Gurjar and D.-Q. Zhang, \emph{{$\pi\sb 1$} of smooth points of a log del
  {P}ezzo surface is finite. {II}}, J. Math. Sci. Univ. Tokyo \textbf{2}
  (1995), no.~1, 165--196. \MR{1348027 (96i:14015)}

\bibitem[K{\etalchar{+}}92]{k-etal}
J.~Koll{\'{a}}r et~al., \emph{Flips and abundance for algebraic threefolds},
  Soc. Math. France, Ast{\'e}risque vol. 211, 1992.

\bibitem[KK10]{k-db}
J{\'a}nos Koll{\'a}r and S{\'a}ndor~J. Kov{\'a}cs, \emph{Log canonical
  singularities are {D}u {B}ois}, J. Amer. Math. Soc. \textbf{23} (2010),
  no.~3, 791--813. \MR{2629988}

\bibitem[KM98]{km-book}
J{\'a}nos Koll{\'a}r and Shigefumi Mori, \emph{Birational geometry of algebraic
  varieties}, Cambridge Tracts in Mathematics, vol. 134, Cambridge University
  Press, Cambridge, 1998, With the collaboration of C. H. Clemens and A. Corti,
  Translated from the 1998 Japanese original. \MR{1658959 (2000b:14018)}

\bibitem[KM99]{kee-mck}
Se{\'a}n Keel and James McKernan, \emph{Rational curves on quasi-projective
  surfaces}, Mem. Amer. Math. Soc. \textbf{140} (1999), no.~669, viii+153.
  \MR{1610249 (99m:14068)}

\bibitem[Kol91]{k-etc}
J{\'a}nos Koll{\'a}r, \emph{Flips, flops, minimal models, etc}, Surveys in
  differential geometry ({C}ambridge, {MA}, 1990), Lehigh Univ., Bethlehem, PA,
  1991, pp.~113--199. \MR{1144527 (93b:14059)}

\bibitem[Kol93]{k-shaf}
\bysame, \emph{Shafarevich maps and plurigenera of algebraic varieties},
  Invent. Math. \textbf{113} (1993), no.~1, 177--215. \MR{1223229 (94m:14018)}

\bibitem[Kol95]{shaf-book}
\bysame, \emph{Shafarevich maps and automorphic forms}, M. B. Porter Lectures,
  Princeton University Press, Princeton, NJ, 1995. \MR{1341589 (96i:14016)}

\bibitem[Mum61]{mumf-top}
David Mumford, \emph{The topology of normal singularities of an algebraic
  surface and a criterion for simplicity}, Inst. Hautes \'Etudes Sci. Publ.
  Math. (1961), no.~9, 5--22. \MR{0153682 (27 \#3643)}

\bibitem[Rei80]{r-c3f}
Miles Reid, \emph{Canonical {$3$}-folds}, Journ\'ees de {G}\'eometrie
  {A}lg\'ebrique d'{A}ngers, {J}uillet 1979/{A}lgebraic {G}eometry, {A}ngers,
  1979, Sijthoff \& Noordhoff, Alphen aan den Rijn, 1980, pp.~273--310.
  \MR{605348 (82i:14025)}

\bibitem[Rei87]{r-ypg}
\bysame, \emph{Young person's guide to canonical singularities}, Algebraic
  geometry, {B}owdoin, 1985 ({B}runswick, {M}aine, 1985), Proc. Sympos. Pure
  Math., vol.~46, Amer. Math. Soc., Providence, RI, 1987, pp.~345--414.
  \MR{927963 (89b:14016)}

\bibitem[Rei94]{reid-nndp}
\bysame, \emph{Nonnormal del {P}ezzo surfaces}, Publ. Res. Inst. Math. Sci.
  \textbf{30} (1994), no.~5, 695--727. \MR{1311389 (96a:14042)}

\bibitem[Sim10]{simp}
Carlos Simpson, \emph{Local systems on proper algebraic {V}-manifolds},
  arXiv1010.3363, 2010.

\bibitem[Tak03]{takayama}
Shigeharu Takayama, \emph{Local simple connectedness of resolutions of
  log-terminal singularities}, Internat. J. Math. \textbf{14} (2003), no.~8,
  825--836. \MR{2013147 (2004m:14023)}

\end{thebibliography}

\newcommand{\etalchar}[1]{$^{#1}$}
\def\cprime{$'$}
\providecommand{\bysame}{\leavevmode\hbox to3em{\hrulefill}\thinspace}
\providecommand{\MR}{\relax\ifhmode\unskip\space\fi MR }
\providecommand{\MRhref}[2]{%
  \href{http://www.ams.org/mathscinet-getitem?mr=#1}{#2}
}
\providecommand{\href}[2]{#2}


\noindent Princeton University, Princeton NJ 08544-1000

{\begin{verbatim}kollar@math.princeton.edu\end{verbatim}}

\end{document}